\documentclass[a4paper,11pt]{article}

\usepackage{moreverb,url}
\usepackage[
	colorlinks,
	bookmarksopen,
	bookmarksnumbered,
	citecolor=blue,
	urlcolor=black,
	linkcolor=blue,
]{hyperref}
\usepackage[font=small,labelfont=bf]{caption}
\usepackage{subcaption}
\usepackage{csquotes}
\usepackage{mathtools}
\usepackage{xcolor}
\usepackage{amsfonts}
\usepackage{amsmath}
\DeclareMathSymbol{\shms}{\mathbin}{AMSa}{"39}  
\usepackage{etoolbox}	
\usepackage{adjustbox}
\usepackage{placeins}
\usepackage[title]{appendix}
\usepackage{float}
\usepackage{siunitx}
\usepackage{breqn}

\usepackage{authblk}		
\usepackage{orcidlink}

\usepackage{geometry}
\usepackage[numbers,sort]{natbib}
\newtheorem{example}{Example}[section]
\newtheorem{remark}{Remark}[section]
\renewcommand{\d}{{\,\rm  d}}

\newcommand\BibTeX{{\rmfamily B\kern-.05em \textsc{i\kern-.025em b}\kern-.08em
			T\kern-.1667em\lower.7ex\hbox{E}\kern-.125emX}}

\usepackage{accents}
\newlength{\dhatheight}

\newcommand{\transp}{^\mathrm{T}}
\newcommand{\ntransp}{^\mathrm{-T}}

\renewcommand{\d}[1][]{\ensuremath{\,\mathrm{d}#1}}
\renewcommand{\phi}{\varphi}
\newcommand{\D}[1][]{\ensuremath{\mathrm{D}#1}}
\usepackage{color}
\usepackage{transparent}
\newcommand{\npe}{^{n+1}}
\newcommand{\npeh}{^{n+1/2}}
\newcommand{\n}{^{n}}

\newcommand{\PDG}[1][]{\ensuremath{\bar{\partial}#1}}
\newcommand{\DG}[1][]{\ensuremath{\bar{\mathrm{D}}#1}}
\newcommand{\eqcomma}{\ensuremath{\text{,}}}
\newcommand{\eqdot}{\ensuremath{\text{.}}}
\renewcommand{\vec}[1]{\ensuremath{\mbox{$ #1 $}}}
\usepackage{bm}
\renewcommand{\epsilon}{\varepsilon}
\newsavebox{\topleftimage}
\newsavebox{\bottomleftimage}
\newsavebox{\bigimage}
\usepackage{mwe}
\usepackage[ruled]{algorithm2e}

\usepackage{pgfplots}
\pgfplotsset{compat = newest}
\usepackage{pgfplotstable}
\usepackage{tikz}
\usepackage{tikz-3dplot} 
\usetikzlibrary{arrows.meta,decorations.markings}
\usepgfplotslibrary{fillbetween}
\pgfkeys{/pgf/number format/.cd,fixed,precision=4}

\newlength\figH
\newlength\figW
\usepgfplotslibrary{patchplots}

\tikzset{myarrowhead/.style={decoration={markings,mark=at position 1 with %
{\arrow[scale=1.,>={Triangle[length=8pt, width=7pt]}]{>}}},postaction={decorate}}}

\DeclareUnicodeCharacter{2212}{−}
\usepgfplotslibrary{groupplots,dateplot}
\usepgfplotslibrary{polar}
\usetikzlibrary{positioning}
\usetikzlibrary{patterns,shapes.arrows}
\pgfplotsset{compat=newest}

\pgfplotsset{every tick label/.append style={font=\scriptsize}}

\definecolor{color1}{RGB}{230, 159, 0}
\definecolor{color2}{RGB}{86, 180, 233}
\definecolor{color3}{RGB}{204, 121, 167}
\definecolor{color4}{RGB}{0, 158, 115}
\definecolor{color5}{RGB}{0, 114, 178}
\definecolor{color6}{RGB}{213, 94, 0}
\definecolor{color7}{RGB}{240, 228, 66}
\definecolor{colorblack}{RGB}{0, 0, 0}          

\definecolor{cgrey}{RGB}{128, 128, 128}
\definecolor{cdarkgrey}{RGB}{98, 98, 98}

\colorlet{corange}{color1}
\colorlet{ccyan}{color2}
\colorlet{cviolet}{color3}
\colorlet{cgreen}{color4}
\colorlet{cblue}{color5}
\colorlet{cred}{color6}
\colorlet{cyellow}{color7}

\colorlet{ctop}{color1}
\colorlet{ccenter}{color2}
\colorlet{cbottom}{color3}

\colorlet{cpath1}{color1}
\colorlet{cpath2}{color2}
\colorlet{cpath3}{color3}
\colorlet{cpath4}{color4}
\colorlet{cpath5}{color7}

\colorlet{cpathneutral}{cgrey}

\colorlet{cpathm90}{color1}
\colorlet{cpathm45}{color2}
\colorlet{cpath0}{color3}
\colorlet{cpath45}{color4}
\colorlet{cpath90}{color7}

\colorlet{cmtoa}{cred}
\colorlet{cHS01}{cgreen}
\colorlet{cHS02}{color4}
\colorlet{cHS03}{cblue}
\colorlet{cmtlinear}{cyellow}
\colorlet{cmtlinearcompliance}{cviolet}

\graphicspath{{./images/}}

\newtoggle{expensiveplots}

\toggletrue{expensiveplots}

\AddToHook{shipout/foreground}{
       \begin{tikzpicture}[overlay, remember picture]
            \node[fill=none, draw, text=red, rotate=0, scale=1, yshift=13.7cm] at (current page) {
				 Note that this article has an open-access published version at \href{https://doi.org/10.1002/pamm.202300144}{\color{red} \textsf{ https://doi.org/10.1002/pamm.202300144}} .};
        \end{tikzpicture}    
    }

\begin{document}

\author[,1]{Philipp L. Kinon \orcidlink{0000-0002-4128-5124}\footnote{Corresponding author: \href{mailto:philipp.kinon@kit.edu}{philipp.kinon@kit.edu}, phone +49 721 608 46081}}
\author[2]{Tobias Thoma \orcidlink{0000-0002-6876-5662}}
\author[1]{Peter Betsch \orcidlink{0000-0002-0596-2503}}
\author[2]{Paul Kotyczka \orcidlink{0000-0002-6669-6368}}

\affil[1]{Institute of Mechanics, Karlsruhe Institute of Technology (KIT), Otto-Ammann-Platz 9, 76131 Karlsruhe, Germany}

\affil[2]{TUM School of Engineering and Design, Technical University of Munich (TUM), Boltzmannstr. 15, 85748 Garching, Germany}
\renewcommand\Authands{ and }

\newcommand{\mypapertitel}{%
	Discrete nonlinear elastodynamics in a port-Hamiltonian framework\footnote{Submitted to Proceedings in Applied Mathematics and Mechanics}
}
\title{\mypapertitel{}}
\maketitle
\begin{abstract}
	\noindent
	We provide a fully nonlinear port-Hamiltonian formulation for discrete elastodynamical systems as well as a structure-preserving time discretization. The governing equations are obtained in a variational manner and represent index-1 differential algebraic equations. Performing an index reduction one obtains the port-Hamiltonian state space model, which features the nonlinear strains as an independent state next to position and velocity. Moreover, hyperelastic material behavior is captured in terms of a nonlinear stored energy function. The model exhibits passivity and losslessness and has an underlying symmetry yielding the conservation of angular momentum. We perform temporal discretization using the midpoint discrete gradient, such that the beneficial properties are inherited by the developed time stepping scheme in a discrete sense. The numerical results obtained in a representative example are demonstrated to validate the findings.
\end{abstract}


\section{Introduction}
Due to its well-known beneficial properties \cite{duindam_modeling_2009}, the port-Hamiltonian (PH) framework has gained large popularity also when modeling elastodynamical systems (which might occur in structural mechanics \cite{warsewa2021port} or multibody system dynamics \cite{brugnoli2021port}). The developed models are able to describe multiphysical coupling as well as the interconnection to other subsystems, whilst achieving a systematic and power-preserving formulation.
However, in the field of elastodynamics quite often the focus has been on linear descriptions (either due to some kinematic assumptions or by considering purely linear-elastic material only).
In the present work, we show a PH formulation for elastodynamics which can capture both geometric and material nonlinearities. While we restrict ourselves to finite-dimensional systems in this contribution, also infinite-dimensional problems fit into this model class after a suitable mixed finite element discretization \cite{thoma2022port,kinon_porthamiltonian_2023}.

This work starts by modeling spatially discrete mechanical systems with hyperelastic elements and introduces a PH model for these problems in Sect.~\ref{sec_model}. Subsequently, a time discretization by means of the midpoint discrete gradient \cite{gonzalez_time_1996} is performed in Sect.~\ref{sec_discrete} which inherits conservation properties of the underlying time-continuous model. Numerical results using the example of a fully nonlinear spring pendulum are shown in Sect.~\ref{sec_numerics}. Section~\ref{sec_conclusion} concludes the findings and gives a brief outlook.

\section{Finite-dimensional nonlinear elastodynamics}\label{sec_model}
In this section, the fundamental governing equations for finite-dimensional (i.e., spatially discrete) elastodynamical systems are outlined. Furthermore, we perform an index reduction of the index-1 differential-algebraic equations in order to provide a new port-Hamiltonian formulation for the problem class at hand.

\subsection{Governing equations}
Consider the motion of a finite-dimensional elastodynamical system consisting of $N$ point masses in three dimensions with the configuration space $\mathcal{Q} \subseteq \mathbb{R}^{3N}$ such that $\vec{q} \in \mathcal{Q}$ denotes the position and $\vec{v} \in T_{\vec{q}} \mathcal{Q} \subseteq \mathbb{R}^{3N}$ the velocity, respectively. Moreover, we introduce strain measures $C \in \mathbb{R}^{n_{\mathrm{el}}}$ corresponding to ${n_{\mathrm{el}}}$ hyperelastic elements in the system and analyze the motion over some time interval $T = [ t_0, t_{\mathrm{end}}] \subset \mathbb{R}_{\geq 0}$. We derive the equations of motion of the system by means of the Lagrange-d'Alembert principle (see, e.g., \cite{marsden_introduction_1999}) such that
\begin{align}
    \delta \int_{t_0}^{t_{\mathrm{end}}} L_{\mathrm{aug}} (\vec{q},\vec{v},\vec{p},\vec{C}, S) \d t +  \int_{t_0}^{t_{\mathrm{end}}} \vec{f}_{\mathrm{np}}(\vec{q},\vec{v}) \cdot \delta \vec{q} \d t = 0 ,\label{livens}
\end{align}
where $\vec{f}_{\mathrm{np}}$ are non-potential forces and the augmented Lagrangian $ L_{\mathrm{aug}} : \mathcal{Q} \times T_{\vec{q}} \mathcal{Q} \times T_{\vec{q}}^* \mathcal{Q} \times \mathbb{R}^{n_{\mathrm{el}}} \times \mathbb{R}^{n_{\mathrm{el}}} \rightarrow \mathbb{R} $ is given by
\begin{align} \label{augm_lagrange}
    L_{\mathrm{aug}} (\vec{q},\vec{v},\vec{p},C, S) = T(\vec{v}) - V(\vec{q},C) + \vec{p} \cdot (\dot{\vec{q}} - \vec{v}) + \frac{1}{2}S \cdot g(\vec{q},C)
\end{align}
Therein, the kinetic (co-)energy $T: T_{\vec{q}} \mathcal{Q} \rightarrow \mathbb{R}$, $\vec{v} \mapsto \frac{1}{2} \vec{v} \cdot \vec{M}\vec{v}$, where $\vec{M} \in \mathbb{R}^{3N \times 3N}$ represents the constant, positive definite and symmetric mass matrix, and the potential energy
$ V(\vec{q},C) = V_{\mathrm{int}}(C) + V_{\mathrm{ext}}(\vec{q})$
are introduced. The potential energy is split up into an internal part $V_{\mathrm{int}}: \mathbb{R}^{n_{\mathrm{el}}} \rightarrow \mathbb{R}$
and an external part $V_{\mathrm{ext}}: \mathcal{Q} \rightarrow \mathbb{R}$.
The variational principle introduces Lagrange multipliers $\vec{p} \in T_{\vec{q}}^* \mathcal{Q} \subseteq \mathbb{R}^{3N}$ in the spirit of Livens principle \cite{livens_hamiltons_1919} in order to enforce the kinematic relation $ \dot{\vec{q}} = \vec{v}$ and stress-like multipliers $S \in \mathbb{R}^{n_{\mathrm{el}}}$ to realize the kinematic side condition
\begin{align} \label{kin_constraint}
    g(\vec{q},\vec{C}) = \vec{C} - \tilde{\vec{C}}(\vec{q}) = 0  \qquad \Leftrightarrow \qquad \begin{bmatrix}
        g_1    \\
        g_2    \\
        \vdots \\
        g_{n_{\mathrm{el}}}
    \end{bmatrix} = \begin{bmatrix}
        C_1    \\
        C_2    \\
        \vdots \\
        C_{n_{\mathrm{el}}}
    \end{bmatrix} - \begin{bmatrix}
        \tilde{C}_1(\vec{q}) \\
        \tilde{C}_2(\vec{q}) \\
        \vdots               \\
        \tilde{C}_{n_{\mathrm{el}}}(\vec{q})
    \end{bmatrix} = \vec{0}
\end{align}
which can be regarded as an operation of column vectors comprising the information of all ${n_{\mathrm{el}}}$ elastic elements. Note that $\tilde{\vec{C}}$ denotes the strain computed by means of the position.
Computing individual stationary conditions pertaining to \eqref{livens} gives rise to the index-1 differential-algebraic equations (DAEs)
\begin{subequations}\label{EL}
    \begin{align}
        \dot{\vec{q}} & = v \label{EL1}  ,                                                                                                                    \\
        \dot{\vec{p}} & = - \D V_{\mathrm{ext}}(\vec{q}) - \frac{1}{2} \D \tilde{C}(\vec{q})\transp S + \vec{f}_{\mathrm{np}}(\vec{q},\vec{v}) ,  \label{EL2} \\
        \vec{p}       & = \vec{M} \vec{v}                                 ,           \label{EL3}                                                             \\
        \frac{1}{2}S  & = \D V_{\mathrm{int}}(C)                         ,    \label{EL4}                                                                     \\
        C             & = \tilde{C}(\vec{q}), \label{EL5}
    \end{align}
\end{subequations}
which are the Euler-Lagrange equations of the problem at hand. Note that $\D (\bullet)$ represents the gradient operator. From relation \eqref{EL3} it becomes visible that the Lagrange multipliers $\vec{p}$ emerge as the conjugate momenta. For further details on Livens principle and an extension to constraint mechanical systems, see \cite{kinon_ggl_2023, kinon_structure_2023}. Furthermore, \eqref{EL4} exposes that $S$ can be viewed as the work-conjugated stresses corresponding to the strains, which are computed as the gradient of the internal energy.

\subsection{Constitutive modeling}\label{sec_2_2}

In analogy to the stored energy density function for geometrically exact strings with hyperelastic material \cite{strohle2022simultaneous,kinon_porthamiltonian_2023} we assume the internal energy as the sum of contributions by all elastic element such that
\begin{align} \label{stored_energy_function}
    V_{\mathrm{int}}(C) = \sum_{i=1}^{{n_{\mathrm{el}}}} V_{\mathrm{int},i}(C_i) = \sum_{i=1}^{{n_{\mathrm{el}}}} \frac{k_i l_{0,i}}{4} \left(C_i - \ln(C_i) - 1 \right) ,
\end{align}
where $k_i$ denotes a spring stiffness parameter and $l_{0,i}$ is the natural length of the $i$-th spring element. It can be verified that the stress-free configuration of the elements is indeed given by $C_i=1$.
Thus, the global minimum of the internal potential is given by $\min(V_{\mathrm{int}}) = \sum_{i=1}^{{n_{\mathrm{el}}}} V_{\mathrm{int},i}(C_i=1) =0$.
Without loss of generality, we consider Cauchy-Green strains
\begin{align}\label{green_lagrange}
    \tilde{C}_i(\vec{q}) = \frac{\bar{\vec{q}}_i \cdot \bar{\vec{q}}_i  }{ l_{0,i}^2} .
\end{align}
where $\bar{\vec{q}}_i$ denotes the relative vector between the two ends of the elastic element. It can be stated that the relative vector between the ends of an elastic element can be obtained from a linear combination of the $N$ positions $\vec{q}_j \in \mathbb{R}^3$  such that $\vec{q}=(\vec{q}_1, \vec{q}_2, \dots, \vec{q}_N)$ and
\begin{align} \label{linear_comb}
    \bar{\vec{q}}_i = \sum_{j=1}^{N} a_{ij} \vec{q}_j
\end{align}
with constant $a_{ik} \in \mathbb{R}$.

\begin{figure}[H]
    \centering
    \vspace*{2.5mm}
    \def\svgwidth{0.4\textwidth}
\begingroup%
  \makeatletter%
  \providecommand\color[2][]{%
    \errmessage{(Inkscape) Color is used for the text in Inkscape, but the package 'color.sty' is not loaded}%
    \renewcommand\color[2][]{}%
  }%
  \providecommand\transparent[1]{%
    \errmessage{(Inkscape) Transparency is used (non-zero) for the text in Inkscape, but the package 'transparent.sty' is not loaded}%
    \renewcommand\transparent[1]{}%
  }%
  \providecommand\rotatebox[2]{#2}%
  \newcommand*\fsize{\dimexpr\f@size pt\relax}%
  \newcommand*\lineheight[1]{\fontsize{\fsize}{#1\fsize}\selectfont}%
  \ifx\svgwidth\undefined%
    \setlength{\unitlength}{242.79505501bp}%
    \ifx\svgscale\undefined%
      \relax%
    \else%
      \setlength{\unitlength}{\unitlength * \real{\svgscale}}%
    \fi%
  \else%
    \setlength{\unitlength}{\svgwidth}%
  \fi%
  \global\let\svgwidth\undefined%
  \global\let\svgscale\undefined%
  \makeatother%
  \begin{picture}(1,0.83576585)%
    \lineheight{1}%
    \setlength\tabcolsep{0pt}%
    \put(0,0){\includegraphics[width=\unitlength,page=1]{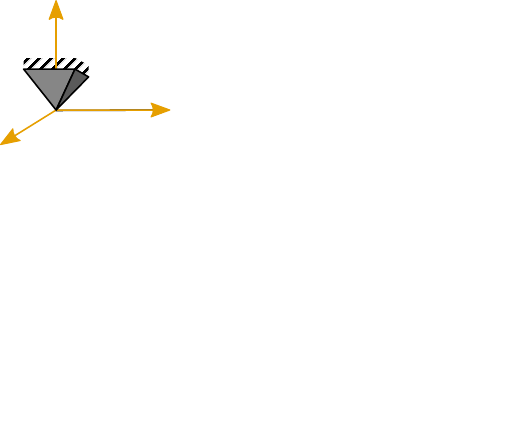}}%
    \put(0.69516637,0.2548996){\color[rgb]{0,0,0}\makebox(0,0)[lt]{\lineheight{1.25}\smash{\begin{tabular}[t]{l}$\vec{q}$\end{tabular}}}}%
    \put(0.46640748,0.47170693){\color[rgb]{0.5254902,0.5254902,0.5254902}\makebox(0,0)[lt]{\lineheight{1.25}\smash{\begin{tabular}[t]{l}$k,l_0$\end{tabular}}}}%
    \put(0,0){\includegraphics[width=\unitlength,page=2]{images/pendulum2.pdf}}%
    \put(0.21557531,0.65968902){\color[rgb]{0.90196078,0.62352941,0}\makebox(0,0)[lt]{\lineheight{1.25}\smash{\begin{tabular}[t]{l}$\{\vec{e}_i\}$\end{tabular}}}}%
    \put(0,0){\includegraphics[width=\unitlength,page=3]{images/pendulum2.pdf}}%
    \put(0.8616892,0.08727617){\color[rgb]{0,0,0}\makebox(0,0)[lt]{\lineheight{1.25}\smash{\begin{tabular}[t]{l}$m$\end{tabular}}}}%
  \end{picture}%
\endgroup%

    \caption{Spring pendulum system with orthonormal basis $\{ \vec{e}_i \}$.}
    \label{fig:sketch}
\end{figure}

\begin{example}\label{ex:pendulum}
    As an example consider the nonlinear spring pendulum with mass $m$ in three dimensions (see, e.g., \cite{betsch_energy_2016}) as depicted in Fig.~\ref{fig:sketch}. Correspondingly, this yields $N=1$ such that $\mathcal{Q} \subseteq \mathbb{R}^3$, $T_{\vec{q}} \mathcal{Q} \subseteq \mathbb{R}^3$ and $T_{\vec{q}}^* \mathcal{Q} \subseteq \mathbb{R}^3$. Moreover, $\vec{M}= m \vec{I} \in \mathbb{R}^{3 \times 3}$. There is one hyperelastic spring element yielding the internal potential governed by \eqref{stored_energy_function} (with  ${n_{\mathrm{el}}}=1$) along with \eqref{stored_energy_function} and $\bar{\vec{q}}_1 = \vec{q}$. The external potential
    \begin{align}
        V_{\mathrm{ext}}(\vec{q}) = - m \vec{b} \cdot \vec{q} \eqdot
    \end{align}
    due to constant gravitational forces $\vec{b}=-9.81 \vec{e}_3$ is considered.
\end{example}

\subsection{Port-Hamiltonian formulation}

The total energy of the system is formulated in its state variables $\vec{x} = (\vec{q}, \vec{v}, \vec{C})$ such that
\begin{align}
    \label{hamiltonian}
    H(\vec{x})=H(\vec{q}, \vec{v},C) = T(\vec{v}) + V_{\mathrm{int}}(C) + V_{\mathrm{ext}}(\vec{q}) .
\end{align}
The co-state variables (efforts) emerge as the partial derivatives
\begin{align}
    \frac{\partial H}{\partial \vec{q}} & = \D V_{\mathrm{ext}}(\vec{q})  \eqcomma   \quad
    \frac{\partial H}{\partial \vec{v}}  = \vec{M} \vec{v} = \vec{p} \eqcomma    \quad
    \frac{\partial H}{\partial C}        = \D V_{\mathrm{int}}(C) = \frac{1}{2}S \eqdot
    \label{efforts}
\end{align}
The strain relation \eqref{EL5} induces the hidden kinematic relation
\begin{align}
    \dot{C} = \D \tilde{C}(\vec{q}) \dot{\vec{q}} = \D \tilde{C}(\vec{q}) \vec{v} \label{kinematic_constraint} .
\end{align}
Next, we perform an index reduction by replacing \eqref{EL5} with its kinematic counterpart \eqref{kinematic_constraint}. Moreover by eliminating the conjugate momenta using \eqref{EL3}, one eventually obtains the state equations governing the finite-dimensonial PH system, which can be written as ordinary differential equations of first order in time
\begin{align}
    \label{eqn_of_motion}
    \begin{bmatrix}
        \vec{I} & 0       & 0       \\
        0       & \vec{M} & 0       \\
        0       & 0       & \vec{I}
    \end{bmatrix}
    \begin{bmatrix}
        \dot{\vec{q}} \\ \dot{\vec{v}} \\ \dot{C}
    \end{bmatrix} & = \begin{bmatrix}
        0        &  & \vec{I}               &  & 0                              \\
        -\vec{I} &  & 0                     &  & - \D \tilde{C}(\vec{q})\transp \\
        0        &  & \D \tilde{C}(\vec{q}) &  & 0
    \end{bmatrix}
    \begin{bmatrix}
        \D V_{\mathrm{ext}}(\vec{q}) \\
        \vec{v}                      \\
        \frac{1}{2}S
    \end{bmatrix} + \begin{bmatrix}
        \vec{0} \\ \vec{I} \\ \vec{0}
    \end{bmatrix} \vec{f}_{\mathrm{np}}(\vec{q},\vec{v}) \eqdot
\end{align}
with $\frac{1}{2}S = \D V_{\mathrm{int}}(C)$. These equations can be recast in the compact form
\begin{subequations}
    \begin{align} \label{eq:PH-continuous-Mehrmann-Morandin2019}
        \vec{E}\dot{\vec{x}} & = \vec{J}(\vec{x}) \vec{z} + \vec{B}\vec{u} \ ,
        \quad\mbox{where}\quad
        \vec{E}\transp\vec{z} = \D H(\vec{x}) \ ,                              \\
        \vec{y}              & = \vec{B}\transp \vec{z},
    \end{align}
\end{subequations}
where $u$ and $y$ denote collocated input and outputs, respectively (cf. \cite{beattie2018linear} for the initial introduction for linear systems, \cite{mehrmann19} for the extension to arbitrary nonlinear Hamiltonians and \cite{brugnoli2021port} for the application in the multibody system context). Note that $\vec{E}\transp = \vec{E}$ and $\vec{J}(\vec{x}) = -\vec{J}(\vec{x})\transp$. While the state-dependence of $\vec{J}$ reflects the geometric nonlinearity, the material nonlinearity due to hyperelastic material behavior is covered in the nonlinear effort law.
\begin{remark}
    Alternatively, a canonical PH representation of the problem at hand can be obtained with a state transformation $\tilde{x}=\vec{E}\vec{x}$, with $E= \D \tilde{x}(x)$, such that the momenta $p$ appear in the state vector $\tilde{\vec{x}}=(\vec{q},\vec{p},\vec{C})$ instead of the velocities. Moreover, the energy function is given by $H = \tilde{H}(\tilde{\vec{x}})$ and $\D \tilde{H}(\tilde{\vec{x}}) = \vec{E}\ntransp \D H (x) = (\D V_{\mathrm{ext}}(\vec{q}), \vec{v}, \frac{1}{2}S)$. Correspondingly, \eqref{eq:PH-continuous-Mehrmann-Morandin2019} would boil down to the canonical representation $\dot{\tilde{x}}= \vec{J}(\tilde{\vec{x}}) \D \tilde{H}(\tilde{\vec{x}}) + \vec{B}\vec{u}$ with the same matrices $\vec{J}$ and $\vec{B}$ as in \eqref{eqn_of_motion}.
\end{remark}
In order to close the PH system representation, initial conditions
\begin{align}
    \vec{q}(t=0) & = \vec{q}_0 , \quad
    \vec{v}(t=0)  = \vec{v}_0  ,\quad
    C(t=0)        = C_0 = \tilde{C}(\vec{q}_0)
\end{align}
are required. We can write the energy balance as
\begin{align}
    \dot{H} & = \D H(\vec{x}) \transp \dot{\vec{x}} = \vec{z} \transp E \dot{\vec{x}} = \vec{z} \transp ( \vec{J}(\vec{x}) \vec{z} + \vec{B}\vec{u}) = \vec{y}\transp \vec{u} ,
    \label{power_balance}
\end{align}
which demonstrates the passivity and losslessness of the system, due to the skew-symmetry of the structure-matrix. In the case of vanishing non-potential forces (i.e., no inputs) the system is energy-preserving such that $\dot{H}=0$. In the absence of system inputs we can moreover analyze the preservation of total angular momentum
\begin{align}
    \vec{L}(\vec{q},\vec{v}) = \sum_{k=1}^{N} \vec{q}_k \times m_k \vec{v}_k
\end{align}
where we decompose positions and corresponding velocities into $N$ three-dimensional contributions such that $\vec{M}_k = m_k \vec{I}_{3 \times 3}$. The balance of linear momentum for this $k$-contribution derived from the second line of \eqref{eqn_of_motion} is given by
\begin{align} \label{particle_momentum_balance}
    m_k \dot{\vec{v}}_k = -\D V^{\mathrm{ext}}_k(\vec{q}_k) - \sum_{i=1}^{{n_{\mathrm{el}}}} \frac{S_i}{2} \D \tilde{\vec{C}}_i(\vec{q}_k)
\end{align}
In view of \eqref{green_lagrange} and \eqref{linear_comb}, one obtains
\begin{align}\label{particle_derivative}
    \D \tilde{\vec{C}}_i(\vec{q}_k) = a_{ik} \frac{2}{l_{0,i}^2} \sum_{j=1}^{N} a_{ij} \vec{q}_j .
\end{align}
Eventually, the time derivative of the total angular momentum can be computed such that
\begin{align}
    \frac{\d}{\d t}\vec{L} & = \sum_{k=1}^{N} ( \dot{\vec{q}}_k \times m_k \vec{v}_k + \vec{q}_k \times m_k \dot{v}_k )  = \sum_{k=1}^{N} ( \vec{v}_k \times m_k \vec{v}_k + \vec{q}_k \times ( - \D V^{\mathrm{ext}}_k(\vec{q}_k) - \sum_{i=1}^{{n_{\mathrm{el}}}} S_i a_{ik} \frac{2}{l_{0,i}^2} \sum_{j=1}^{N} a_{ij} \vec{q}_j ) )                                     \notag               \\
                           & = - \sum_{k=1}^{N} \vec{q}_k \times \D V^{\mathrm{ext}}_k(\vec{q}_k) - \sum_{i=1}^{{n_{\mathrm{el}}}} \frac{2 S_i}{l_{0,i}^2} \sum_{k=1}^{N} a_{ik} \vec{q}_k \times \sum_{j=1}^{N} a_{ij} \vec{q}_j                                                                                          = - \sum_{k=1}^{N} \vec{q}_k \times \D V^{\mathrm{ext}}_k(\vec{q}_k)
\end{align}
where \eqref{particle_momentum_balance}-\eqref{particle_derivative} have been taken into account along with $\dot{\vec{q}}_k = \vec{v}_k$. Consequently, the angular momentum about the axis of constant external potential forces $\vec{b} = - \D V_{\mathrm{ext}}(\vec{q}_k)$ is preserved, such that
\begin{align}
    \frac{\d}{\d t} \left( \vec{L \cdot \vec{b}} \right) = 0 . \label{ang_mom_cons}
\end{align}
\section{Structure-preserving time discretization}
\label{sec_discrete}

We aim at a structure-preserving time discretization of the PH system \eqref{eqn_of_motion}. To this end, the time interval is discretized with $n_{\mathrm{t}}$ time steps of constant size $h = t\npe-t\n$ such that
$ T = [t_0, t_{\mathrm{end}}] = \bigcup_{n=0}^{n_{\mathrm{t}}} [t\n, t\npe]$
and we apply a one-step scheme which is closely related to the implicit midpoint rule. Let $\vec{x}\n$ be the approximation of the state $\vec{x}(t\n)$ at time $t\n$. The time stepping scheme can now be written in the form
\begin{subequations}
    \label{timediscrete_EoM}
    \begin{equation} \label{timediscrete_EoM1}
        \begin{array}{rcl}
            \vec{E}\left({\vec{x}}\npe-{\vec{x}}\n\right) & = & h ( \vec{J}({\vec{x}}\npeh) {{z}}\npeh + \vec{B} \vec{u}\npeh )  \eqcomma
        \end{array}
    \end{equation}
    where ${\vec{x}}\npeh=\frac{1}{2}({\vec{x}}\npe+{\vec{x}}\n)$ and the collocated outputs in discrete time are
    \begin{align}
        \vec{y}\npeh = \vec{B}\transp \vec{z}\npeh .
    \end{align}
    Moreover, ${\vec{z}}\npeh$ is defined through
    \begin{equation} \label{discrete_gradient_linking}
        \vec{E}\transp{\vec{z}}\npeh = \DG {H}({\vec{x}}\n,{\vec{x}}\npe) \eqcomma
    \end{equation}
\end{subequations}
where $\DG {H}({\vec{x}}\n,{\vec{x}}\npe)$ is the midpoint discrete gradient in the sense of Gonzalez \cite{gonzalez_time_1996}. Among other properties (see \cite{hairer_geometric_2006,gonzalez_time_1996} for more details), $\DG {H}({\vec{x}}\n,{\vec{x}}\npe)$ satisfies the crucial directionality property
\begin{equation} \label{eq:directionality}
    \DG {H}({\vec{x}}\n,{\vec{x}}\npe)\transp \left({\vec{x}}\npe-{\vec{x}}\n\right) = H\npe - H\n ,
\end{equation}
where we have introduced the discrete energy function
$ H\n = H(\vec{x}\n) = T(\vec{v}\n) + V(\vec{q}\n, C\n) = T\n + V\n$.
In particular, for the discrete derivatives we consider
\begin{equation} \label{discrete_derivative}
    \DG {H}({\vec{x}}\n,{\vec{x}}\npe) =
    \begin{bmatrix}
        \PDG_{{\vec{r}}}{H} \\
        \PDG_{{\vec{v}}}{H} \\
        \PDG_{{C}}{H}
    \end{bmatrix} =
    \begin{bmatrix}
        \DG V_{\mathrm{ext}}(\vec{q}\n, \vec{q}\npe) \\
        \vec{M}{\vec{v}}\npeh                        \\
        \DG V_{\mathrm{int}}(C\n, C\npe)
    \end{bmatrix} \eqcomma
\end{equation}
where the classical Greenspan's formula \cite{greenspan_conservative_1984} can be applied componentwise in the third block entry. In the limit $C\npe \rightarrow C\n$, the midpoint evaluation of the analytical gradient is used.
It can be easily checked by a straightforward calculation that \eqref{discrete_derivative} does indeed satisfy the directionality property \eqref{eq:directionality}.
We can now analyze the scheme with respect to conservation principles. Firstly, making use of \eqref{timediscrete_EoM} and directionality property \eqref{eq:directionality} we obtain
\begin{align} \label{discrete_energy_cons}
    {H}\npe - {H}\n \notag
    =\vec{z}\npeh \cdot \vec{E}\left({\vec{x}}\npe-{\vec{x}}\n\right)
    ={\vec{z}}\npeh \cdot h  \vec{B}\vec{u}\npeh
    = h\, \vec{y}\npeh \cdot \vec{u}\npeh \eqcomma
\end{align}
which is a discrete counterpart of \eqref{power_balance}. This proves that the present time-stepping scheme exhibits passivity, losslessness and for $u=0$ also energy-conservation. Note that these properties are a direct consequence of the skew symmetry of $\vec{J}(\vec{x})$ and hold regardless of the evaluation point of $\vec{J}(\vec{x})$ in time. However, a midpoint evaluation is crucial for the scheme to be symmetric which is said to enhance robustness (see \cite{hairer_geometric_2006}). Moreover, the conservation of further discrete quantities relies on this discretization choice. For example, consider the discrete, total angular momentum
\begin{align}
    \vec{L}\n = \sum_{k=1}^{N} \vec{q}\n_k \times m_k \vec{v}\n_k \eqdot
\end{align}
Using the first and second equation in \eqref{timediscrete_EoM} one obtains
\begin{align}
    \vec{L}\npe - \vec{L}\n & =
    \sum_{k=1}^{N} \left( (\vec{q}\npe_k - \vec{q}\n_k) \times m_k \vec{v}\npeh_k + \vec{q}\npeh_k \times m_k (\vec{v}\npe_k - \vec{v}\n_k) \right) \notag \\  &=
    - h \sum_{k=1}^{N} \left( \DG V^{\mathrm{ext}}_k(\vec{q}\n_k, \vec{q}\npe_k) + \sum_{i=1}^{{n_{\mathrm{el}}}} \D \tilde{C}_i(\vec{q}\npeh_k)\transp \frac{1}{2}\bar{S}_i  \right) \times \vec{q}\npeh_k ,
\end{align}
where $ \frac{1}{2}\bar{S}_i = \DG V^{\mathrm{int}}_i(C_i\n, C_i\npe)$ and colinear velocites have been taken into account. Moreover, the second term on the right hand side vanishes due to \eqref{particle_derivative} evaluated in the midpoint. Note that the midpoint evaluation of the skew-symmetric matrix $J$ in \eqref{timediscrete_EoM} is crucial here. Eventually, the angular momentum component about the axis of external potential forces $\bar{\vec{b}} = - \DG V^{\mathrm{ext}}_k(\vec{q}\n_k, \vec{q}\npe_k)$ for all $k$ is a conserved quantity also in discrete time, such that
\begin{align}
    (\vec{L}\npe - \vec{L}\n) \cdot \bar{\vec{b}} & = 0 \label{discrete_angmom_cons} .
\end{align}
\begin{remark}
    The present method \eqref{timediscrete_EoM} also covers the kinematic condition \eqref{kin_constraint} exactly on the discrete level, i.e.,
    \begin{align} \label{cons_kin_relation}
        g(\vec{q}\n, C \n) = C\n - \tilde{C}(\vec{q}\n) = 0
    \end{align}
    for all $n=1,\dots,n_{\mathrm{t}}$, if starting with consistent initial conditions satisfying $C_0 = \tilde{C}(\vec{q}_0)$. This property can be shown using the first and third equation in \eqref{timediscrete_EoM1}
    and the fact that $\tilde{C}$ is a quadratic function of its argument in each component.
\end{remark}
\begin{figure}[b]
    \begin{minipage}{0.525\textwidth}
        \centering
        \setlength{\figH}{0.25\textheight}
        \setlength{\figW}{0.2833\textheight}
        \input{images/final_conf.tikz}
        \caption{Initial and final configuration with trajectory during the motion}
        \label{fig:snapshots}
    \end{minipage}
    \hfill
    \begin{minipage}{0.475\textwidth}
        \centering
        \setlength{\figH}{0.15\textheight}
        \setlength{\figW}{0.85\textwidth}
        \input{images/C_diff.tikz}
        \caption{Discrete kinematic condition \eqref{cons_kin_relation}}
        \label{fig:Constitutive_diff}
    \end{minipage}
\end{figure}
\begin{table}[b]
    \centering
    \captionsetup{margin=4cc}
    \caption{Spring pendulum example: Simulation parameters.} \label{tab:tab1}
    \def\arraystretch{1.3}
    \begin{tabular}{@{}lllllllll@{}}
        \hline
        $l_0 \, [\mathrm{m}] $ & $k \, [\mathrm{N}] $ & $ m \, [\mathrm{kg}] $ & $h \, [\mathrm{s}] $ & $T \, [\mathrm{s}]$ & $\epsilon_{\mathrm{Newton}} \, [\mathrm{-}] $ & $\vec{q}_0 \, [\mathrm{m}] $ & $\vec{v}_0 \, [\mathrm{m/s}] $ & $C_0 \,  [\mathrm{-}]     $   \\
        \hline
        $1$                    & $10^4$               & $1$                    & $ 10^{-2}$           & $[0, 4]$            & $10^{-9}$                                     & $(1.1, 0, 0)$                & $(0, 1, 1)$                    & $\tilde{C}(\vec{q}_0) = 1.21$ \\
        \hline
    \end{tabular}
\end{table}

\section{Numerical results}
\label{sec_numerics}

We investigate the motion of a spring pendulum as introduced in Example \ref{ex:pendulum} governed by the equations of motion in PH form \eqref{eqn_of_motion}. We use the time-stepping method given by \eqref{timediscrete_EoM}, which has been solved using
Newton's method down to a residual norm error of $\epsilon_{\mathrm{Newton}}$ in each timestep. The chosen physical and numerical parameters as well as initial conditions are comprised in Table \ref{tab:tab1} and the computer implementation can be found at \cite{kinon_metis_2023}.

Initial and final configuration are displayed in Fig.~\ref{fig:snapshots} together with the computed trajectory of the mass. Moreover, Fig.~\ref{fig:Constitutive_diff} highlights the numerically exact satisfaction of the kinematic condition \eqref{cons_kin_relation}. During the motion energy is heavily exchanged between kinetic and potential parts leading to the evolution, which can be observed in Fig.~\ref{fig:energy}. As proven in \eqref{discrete_energy_cons}, the total energy is preserved numerically, which is verified by means of time increments on the level of computer precision (see Fig.~\ref{fig:energy_diff}). Due to the chosen external potential forces, the angular momentum about the vertical axis is also a conserved quantity (cf. relation \eqref{ang_mom_cons}). The corresponding angular momentum conservation in discrete time \eqref{discrete_angmom_cons} is verified in Fig.~\ref{fig:ang_mom} together with the incremental changes displayed in Fig.~\ref{fig:ang_mom_diff}.

\begin{figure}[htb]
    \begin{minipage}{0.475\textwidth}
        \centering
        \vspace*{3.6mm}
        \setlength{\figH}{0.15\textheight}
        \setlength{\figW}{0.85\textwidth}
        \input{images/energy.tikz}
        \caption{Energy quantities}
        \label{fig:energy}
    \end{minipage}
    \hfill
    \begin{minipage}{0.475\textwidth}
        \centering
        \setlength{\figH}{0.15\textheight}
        \setlength{\figW}{0.85\textwidth}
        \input{images/H_diff.tikz}
        \caption{Increments $H\npe - H\n$}
        \label{fig:energy_diff}
    \end{minipage}
\end{figure}

\begin{figure}[htb]
    \begin{minipage}{0.475\textwidth}
        \centering
        \vspace*{3.6mm}
        \setlength{\figH}{0.15\textheight}
        \setlength{\figW}{0.85\textwidth}
        \input{images/ang_mom.tikz}
        \caption{Angular momentum components}
        \label{fig:ang_mom}
    \end{minipage}
    \hfill
    \begin{minipage}{0.475\textwidth}
        \centering
        \setlength{\figH}{0.15\textheight}
        \setlength{\figW}{0.85\textwidth}
%

%
\begin{tikzpicture}

  \begin{axis}[%
      width=0.951\figW,
      height=\figH,
      at={(0\figW,0\figH)},
      scale only axis,
      xmin=0,
      xmax=4,
      xlabel style={font=\color{white!15!black}},
      xlabel={$t \, [\mathrm{s}]$},
      ymin=-2e-14,
      ymax=2e-14,
      ylabel style={font=\color{white!15!black}},
      ylabel={$\mathrm{[J\,s]}$},
      axis background/.style={fill=white},
      legend style={legend cell align=left, align=left, draw=white!15!black}
    ]

    \addplot [color=color3, line width=1pt]
    table[row sep=crcr]{%
        0	0\\
        0.01	0\\
        0.02	0\\
        0.03	2.22e-16\\
        0.04	-2.22e-16\\
        0.05	0\\
        0.06	-2.22e-16\\
        0.07	0\\
        0.08	0\\
        0.09	0\\
        0.1	-2.22e-16\\
        0.11	4.441e-16\\
        0.12	-2.22e-16\\
        0.13	-4.4411e-16\\
        0.14	-2.22e-16\\
        0.15	4.441e-16\\
        0.16	-6.661e-16\\
        0.17	2.221e-16\\
        0.18	0\\
        0.19	4.441e-16\\
        0.2	-2.22e-16\\
        0.21	-2.22e-16\\
        0.22	0\\
        0.23	4.441e-16\\
        0.24	0\\
        0.25	0\\
        0.26	0\\
        0.27	0\\
        0.28	-4.441e-16\\
        0.29	4.441e-16\\
        0.3	0\\
        0.31	4.441e-16\\
        0.32	4.441e-16\\
        0.33	0\\
        0.34	-4.441e-16\\
        0.35	0\\
        0.36	2.22e-16\\
        0.37	-2.22e-16\\
        0.38	8.882e-16\\
        0.39	-4.441e-16\\
        0.4	2.22e-16\\
        0.41	-2.22e-16\\
        0.42	0\\
        0.43	-8.882e-16\\
        0.44	-4.441e-16\\
        0.45	4.441e-16\\
        0.46	0\\
        0.47	8.882e-16\\
        0.48	0\\
        0.49	-6.661e-16\\
        0.5	0\\
        0.51	-2.22e-16\\
        0.52	4.441e-16\\
        0.53	-2.22e-16\\
        0.54	4.441e-16\\
        0.55	2.22e-16\\
        0.56	0\\
        0.57	-2.22e-16\\
        0.58	0\\
        0.59	2.22e-16\\
        0.6	0\\
        0.61	-4.441e-16\\
        0.62	0\\
        0.63	0\\
        0.64	0\\
        0.65	2.22e-16\\
        0.66	-2.22e-16\\
        0.67	2.22e-16\\
        0.68	0\\
        0.69	0\\
        0.7	0\\
        0.71	2.22e-16\\
        0.72	0\\
        0.73	0\\
        0.74	0\\
        0.75	0\\
        0.76	2.22e-16\\
        0.77	-2.22e-16\\
        0.78	-2.22e-16\\
        0.79	0\\
        0.8	-2.22e-16\\
        0.81	0\\
        0.82	2.22e-16\\
        0.83	0\\
        0.84	0\\
        0.85	2.22e-16\\
        0.86	-2.22e-16\\
        0.87	-2.22e-16\\
        0.88	-2.22e-16\\
        0.89	0\\
        0.9	0\\
        0.91	0\\
        0.92	-6.661e-16\\
        0.93	-2.22e-16\\
        0.94	2.22e-16\\
        0.95	6.661e-16\\
        0.96	0\\
        0.97	6.661e-16\\
        0.98	0\\
        0.99	4.441e-16\\
        1	-1.11e-15\\
        1.01	6.661e-16\\
        1.02	0\\
        1.03	-4.441e-16\\
        1.04	-2.22e-16\\
        1.05	2.22e-16\\
        1.06	-8.882e-16\\
        1.07	2.22e-16\\
        1.08	2.22e-16\\
        1.09	4.441e-16\\
        1.1	-4.441e-16\\
        1.11	-2.22e-16\\
        1.12	-6.661e-16\\
        1.13	0\\
        1.14	0\\
        1.15	1.332e-15\\
        1.16	0\\
        1.17	-4.441e-16\\
        1.18	0\\
        1.19	-8.882e-16\\
        1.2	8.882e-16\\
        1.21	-6.661e-16\\
        1.22	2.22e-16\\
        1.23	8.882e-16\\
        1.24	-4.441e-16\\
        1.25	0\\
        1.26	-8.882e-16\\
        1.27	4.441e-16\\
        1.28	2.22e-16\\
        1.29	-6.661e-16\\
        1.3	0\\
        1.31	-6.661e-16\\
        1.32	-2.22e-16\\
        1.33	-8.882e-16\\
        1.34	6.661e-16\\
        1.35	6.661e-16\\
        1.36	-4.441e-16\\
        1.37	0\\
        1.38	-4.441e-16\\
        1.39	-1.332e-15\\
        1.4	-4.441e-16\\
        1.41	2.22e-16\\
        1.42	-6.661e-16\\
        1.43	8.882e-16\\
        1.44	4.441e-16\\
        1.45	4.441e-16\\
        1.46	0\\
        1.47	4.441e-16\\
        1.48	4.441e-16\\
        1.49	-4.441e-16\\
        1.5	-4.441e-16\\
        1.51	0\\
        1.52	-8.882e-16\\
        1.53	4.441e-16\\
        1.54	6.661e-16\\
        1.55	2.22e-16\\
        1.56	0\\
        1.57	4.441e-16\\
        1.58	-4.441e-16\\
        1.59	-2.22e-16\\
        1.6	2.22e-16\\
        1.61	-4.441e-16\\
        1.62	2.22e-16\\
        1.63	2.22e-16\\
        1.64	0\\
        1.65	0\\
        1.66	-2.22e-16\\
        1.67	0\\
        1.68	4.441e-16\\
        1.69	-2.22e-16\\
        1.7	-2.22e-16\\
        1.71	4.441e-16\\
        1.72	2.22e-16\\
        1.73	0\\
        1.74	0\\
        1.75	0\\
        1.76	0\\
        1.77	0\\
        1.78	0\\
        1.79	0\\
        1.8	-2.22e-16\\
        1.81	0\\
        1.82	2.22e-16\\
        1.83	-2.22e-16\\
        1.84	0\\
        1.85	0\\
        1.86	2.22e-16\\
        1.87	-2.22e-16\\
        1.88	0\\
        1.89	0\\
        1.9	0\\
        1.91	-2.22e-16\\
        1.92	2.22e-16\\
        1.93	2.22e-16\\
        1.94	-2.22e-16\\
        1.95	0\\
        1.96	0\\
        1.97	-2.22e-16\\
        1.98	2.22e-16\\
        1.99	-4.441e-16\\
        2	2.22e-16\\
        2.01	-4.441e-16\\
        2.02	0\\
        2.03	-4.441e-16\\
        2.04	-2.22e-16\\
        2.05	-2.22e-16\\
        2.06	2.22e-16\\
        2.07	2.22e-16\\
        2.08	0\\
        2.09	-4.441e-16\\
        2.1	4.441e-16\\
        2.11	-4.441e-16\\
        2.12	-4.441e-16\\
        2.13	-2.22e-16\\
        2.14	6.661e-16\\
        2.15	0\\
        2.16	0\\
        2.17	-4.441e-16\\
        2.18	0\\
        2.19	2.22e-16\\
        2.2	-6.661e-16\\
        2.21	0\\
        2.22	4.441e-16\\
        2.23	0\\
        2.24	8.882e-16\\
        2.25	0\\
        2.26	-4.441e-16\\
        2.27	6.661e-16\\
        2.28	-4.441e-16\\
        2.29	2.22e-16\\
        2.3	0\\
        2.31	-4.441e-16\\
        2.32	2.22e-16\\
        2.33	2.22e-16\\
        2.34	-2.22e-16\\
        2.35	0\\
        2.36	2.22e-16\\
        2.37	0\\
        2.38	2.22e-16\\
        2.39	0\\
        2.4	0\\
        2.41	-2.22e-16\\
        2.42	0\\
        2.43	0\\
        2.44	-4.441e-16\\
        2.45	4.441e-16\\
        2.46	2.22e-16\\
        2.47	-2.22e-16\\
        2.48	0\\
        2.49	2.22e-16\\
        2.5	-2.22e-16\\
        2.51	2.22e-16\\
        2.52	0\\
        2.53	-2.22e-16\\
        2.54	0\\
        2.55	0\\
        2.56	2.22e-16\\
        2.57	0\\
        2.58	0\\
        2.59	-2.22e-16\\
        2.6	0\\
        2.61	0\\
        2.62	0\\
        2.63	0\\
        2.64	0\\
        2.65	0\\
        2.66	0\\
        2.67	-2.22e-16\\
        2.68	0\\
        2.69	-2.22e-16\\
        2.7	0\\
        2.71	0\\
        2.72	-2.22e-16\\
        2.73	2.22e-16\\
        2.74	0\\
        2.75	2.22e-16\\
        2.76	2.22e-16\\
        2.77	0\\
        2.78	0\\
        2.79	0\\
        2.8	-2.22e-16\\
        2.81	2.22e-16\\
        2.82	0\\
        2.83	0\\
        2.84	-4.441e-16\\
        2.85	4.441e-16\\
        2.86	-4.441e-16\\
        2.87	2.22e-16\\
        2.88	2.22e-16\\
        2.89	0\\
        2.9	-4.441e-16\\
        2.91	4.441e-16\\
        2.92	-6.661e-16\\
        2.93	6.661e-16\\
        2.94	2.22e-16\\
        2.95	-2.22e-16\\
        2.96	0\\
        2.97	0\\
        2.98	0\\
        2.99	0\\
        3	0\\
        3.01	0\\
        3.02	0\\
        3.03	0\\
        3.04	2.22e-16\\
        3.05	-4.441e-16\\
        3.06	2.22e-16\\
        3.07	0\\
        3.08	0\\
        3.09	0\\
        3.1	-2.22e-16\\
        3.11	0\\
        3.12	-2.22e-16\\
        3.13	0\\
        3.14	0\\
        3.15	0\\
        3.16	0\\
        3.17	0\\
        3.18	-2.22e-16\\
        3.19	2.22e-16\\
        3.2	-2.22e-16\\
        3.21	2.22e-16\\
        3.22	0\\
        3.23	0\\
        3.24	2.22e-16\\
        3.25	2.22e-16\\
        3.26	0\\
        3.27	-2.22e-16\\
        3.28	0\\
        3.29	0\\
        3.3	-2.22e-16\\
        3.31	0\\
        3.32	0\\
        3.33	2.22e-16\\
        3.34	2.22e-16\\
        3.35	0\\
        3.36	0\\
        3.37	0\\
        3.38	1.11e-15\\
        3.39	-4.441e-16\\
        3.4	2.22e-16\\
        3.41	0\\
        3.42	-4.441e-16\\
        3.43	0\\
        3.44	2.22e-16\\
        3.45	0\\
        3.46	2.22e-16\\
        3.47	0\\
        3.48	6.661e-16\\
        3.49	-2.22e-16\\
        3.5	-4.441e-16\\
        3.51	-4.441e-16\\
        3.52	-4.441e-16\\
        3.53	0\\
        3.54	-4.441e-16\\
        3.55	8.882e-16\\
        3.56	4.441e-16\\
        3.57	2.22e-16\\
        3.58	6.661e-16\\
        3.59	4.441e-16\\
        3.6	4.441e-16\\
        3.61	-4.441e-16\\
        3.62	-1.332e-15\\
        3.63	4.441e-16\\
        3.64	-2.22e-16\\
        3.65	4.441e-16\\
        3.66	-1.11e-15\\
        3.67	0\\
        3.68	-2.22e-16\\
        3.69	2.22e-16\\
        3.7	4.441e-16\\
        3.71	0\\
        3.72	-4.441e-16\\
        3.73	0\\
        3.74	4.441e-16\\
        3.75	6.661e-16\\
        3.76	-2.22e-16\\
        3.77	4.441e-16\\
        3.78	4.441e-16\\
        3.79	0\\
        3.8	8.882e-16\\
        3.81	0\\
        3.82	-1.332e-15\\
        3.83	1.332e-15\\
        3.84	-4.441e-16\\
        3.85	6.661e-16\\
        3.86	-1.11e-15\\
        3.87	-4.441e-16\\
        3.88	-6.661e-16\\
        3.89	2.22e-16\\
        3.9	-4.441e-16\\
        3.91	0\\
        3.92	4.441e-16\\
        3.93	-8.882e-16\\
        3.94	4.441e-16\\
        3.95	6.661e-16\\
        3.96	-2.22e-16\\
        3.97	0\\
        3.98	6.661e-16\\
        3.99	2.22e-16\\
      };

    \addplot [color=black, dashed, line width=1pt, forget plot]
    table[row sep=crcr]{%
        0	1e-14\\
        4	1e-14\\
      };
    \addplot [color=black, dashed, line width=1pt, forget plot]
    table[row sep=crcr]{%
        0	-1e-14\\
        4	-1e-14\\
      };
  \end{axis}
\end{tikzpicture}%
        \caption{Increments $(L_3)\npe - (L_3)\n$}
        \label{fig:ang_mom_diff}
    \end{minipage}
\end{figure}

\section{Conclusion}
\label{sec_conclusion}
We have presented a PH formulation for the analysis of nonlinear elastodynamical systems with an independent strain quantity in the state vector. This formulation features a non-quadratic energy function in order to account for hyperelastic material behavior. The geometric nonlinearities are completely reflected in the state-dependence of the skew-symmetric matrix operator, while hyperelastic material behavior is contained in the nonlinear constitutive law. We have also shown the passivity of the system (which includes losslessness and energy-conservation). The structure-preserving temporal discretization using midpoint discrete gradients leads to a simulation scheme with discrete time passivity (including losslessness and energy-conservation) and preservation of the angular momentum map.

In the future, the model should be extended to viscous damping and systems with holonomic constraints. Moreover, the feedforward and feedback control design of this model can offer further challenges. Interesting questions might also arise when investigating the model order reduction of this type of nonlinear systems.

\FloatBarrier


\newenvironment{authcontrib}[1]{%
	\subsection*{\textnormal{\textbf{Author Contributions}}}%
	\noindent #1}%
{}%
\newenvironment{acks}[1]{%
	\subsection*{\textnormal{\textbf{Acknowledgements}}}%
	\noindent #1}%
{}%
\newenvironment{funding}[1]{%
	\subsection*{\textnormal{\textbf{Funding}}}%
	\noindent #1}%
{}%
\newenvironment{dci}[1]{%
	\subsection*{\textnormal{\textbf{Declaration of conflicting interests}}}%
	\noindent #1}%
{}%
\newenvironment{code}[1]{%
	\subsection*{\textnormal{\textbf{Code}}}%
	\noindent #1}%
{}%

\begin{acks}
	The financial support for this work by the DFG (German Research Foundation) – project number 388118188 - is gratefully acknowledged.
\end{acks}
\begin{dci}
	The authors declare no conflict of interest.

\end{dci}
\begin{code}
	The source code used for the computations is openly available at \cite{kinon_metis_2023}.
\end{code}

\addcontentsline{toc}{section}{References}
\bibliographystyle{unsrtnat}
\bibliography{bib}

\end{document}